
\input amstex
\documentstyle {amsppt}

\def\S{\frak S}                 
\def\Sym{Sym}            
\def\ss{s_{\heartsuit}}
\def\vv{\upsilon}
\def\d{\partial}
\def\dd{\partial_{\heartsuit}}
\def\l{\lambda} 
\def\t{\tau} 
\def\Pol{{\Cal P}ol}
\def\X{{\Cal X}}
\def\Y{{\Cal Y}}
\def\Z{\Bbb Z}           
\def\N{\Bbb N}           
\def\Q{\widetilde{Q}}
\def\P{\widetilde{P}}
\def\ov{\overline }

\def\semi{\ltimes} 

\def\sv{\sevenrm }
\def\s{\scriptstyle }
\catcode`@=11
\def\m@th{\mathsurround=\z@}
\def\smatrix#1{\null\,\vcenter{\baselineskip=0pt\m@th    
\ialign{\hfil$\s ##$\hfil&&\ \hfil$\s ##$\hfil\crcr      
\mathstrut\crcr\noalign{\kern-\baselineskip}     
#1\crcr\mathstrut\crcr\noalign{\kern-\baselineskip}}}\,}

\magnification 1200

\document
\topmatter
\title{Orthogonal divided differences and Schubert polynomials, 
$\P$-functions, and vertex operators}
\endtitle
\rightheadtext{Lascoux and Pragacz}
\leftheadtext{Orthogonal divided differences and vertex operators}
\author
Alain Lascoux \ and \ Piotr Pragacz
\endauthor
\thanks
This research was supported by grant No. 5031 of French-Polish 
cooperation C.N.R.S.--P.A.N. 
and KBN grant No. 2P03A 05112.
\endthanks
\endtopmatter

\centerline{\it Dedicated to Bill Fulton on his 60th birthday}
\vskip20pt

\leftline{\bf Abstract}
\noindent
{\sv Using divided differences associated with the orthogonal groups,
we investigate the structure of the polynomial rings over
the rings of invariants of the corresponding Weyl groups.
We study in more detail the action of orthogonal divided differences
on some distinguished symmetric polynomials ($\P$-polynomials) and relate
it to vertex operators.
Relevant families of orthogonal Schubert polynomials,
generalizing $\P$-polynomials, and well-suited to intersection
theory computations, are also studied.} 

\vskip20pt

\leftline{\bf Contents:}

\noindent
Introduction. 
\ \ 1.Divided differences. 
\ \ 2.Bases of polynomial rings.   
\ \ 3.Vertex operators. 
\ \ 4.Applications to $\P$-polynomials and orthogonal Schubert polynomials.  
\ \ Appendix: results in type B. 
\ \ References.

\vskip30pt

\head
{\bf Introduction}
\endhead

Divided differences were introduced by Newton 
in his famous interpolation formula (cf. [N, pp.481--483],
and [L] for some historical comments).

Their importance in geometry was shown by Bernstein-Gel'fand-Gel'fand [BGG]
and Demazure [D1,D2] in the context of {\it Schubert calculus} for generalized
flag varieties associated with semisimple algebraic groups in the early 1970's.
More recently, simple divided differences, interpreted as 
{\it correspondences} in flag bundles,
were extensively used in the sequence of papers [F1,F2,F3] by Fulton 
in the context of degeneracy loci associated with classical 
groups. Still another interpretation of divided differences, as 
{\it Gysin morphisms} in the cohomology of flag bundles 
associated with semisimple algebraic groups,
was discussed in [P2, Sect.4] and [PR, Sect.5].
We refer to the lecture notes [FP] for an introduction.  

The case of $SL(n)$ has been developed by the first author 
and Sch\"utzenberger (see e.g. \ [LS1,LS2,LS3] and [M2]).

For other classical groups, parallel studies were done by Billey-Haiman [BH],
Fomin-Kirillov [FK], Ratajski and the second author [PR], and by the authors
[LP1]. The present paper is a continuation of [LP1]. 
Here we study divided differences
associated with the orthogonal groups $SO(2n)$ and $SO(2n+1)$
(i.e. for types D and B).
The results for type B are an immediate adaptation of the results
for type C given in [LP1]. We summarize them in the appendix. We add, however,
a certain new result (Theorem 9) for type B, whose C-analogue was 
not needed in our former paper [LP1].

Our results for type D require some new computations with vertex
operators, that are furnished in Section 3 of the present paper
and summarized in our main Theorem 11. 
In order to simplify the computations with divided differences,
we display them as planar arrays, which allows us to perform
some kind of `` jeu de taquin". This offers a certain technical novelty w.r.t.
our former paper [LP1]. 
In type $C_n$ (or $B_n$), the key r\^ole
was/is played by the divided differences of the form
$$
(\d_0 \d_1 \cdots \d_{n-1}) \cdots (\d_0 \d_1 \cdots \d_{n-k})\,,
\eqno(*)
$$
where $k\le n$.
It appears that in type $D_n$ a similar r\^ole is played by 
the divided differences of the form
$$
(\dd \d_2\cdots \d_{n-1} \d_1 \d_2 \cdots \d_{n-2})
             \cdots  (\dd \d_2\cdots \d_{n-2k+1} \d_1 \d_2 \cdots \d_{n-2k})\,,
\eqno(**)
$$
where $k\le n/2$. \
Here $\d_i$, for $i>0$, are Newton's (simple) divided differences:
$$
f\, \d_i :=  {{f-f(\ldots,x_{i+1},x_i,\ldots)}\over{x_i-x_{i+1}}}\,,
$$
and moreover we set
$$ 
f\, \d_0 :=  {{f-f(-x_1,x_2,\ldots)}\over{-x_1}}\,,   
$$
$$ 
f\, \dd := {{f-f(-x_2,-x_1,x_3,\ldots)}\over{-x_1 - x_2}}\,.
$$
We compose the simple orthogonal divided differences in (*) and (**) 
{\it from left to right}. As the the Weyl group of type $D$
is naturally embedded in the Weyl group of type $B$, the divided
difference (**) can be expressed in terms of (*). Such basic relations
are given in Proposition 6 and Corollary 8. 

The symmetric functions which are most adapted to orthogonal divided
differences are $\P$-polynomials [PR], which are a variant of
Schur $P$-polynomials.

Our paper is of an algebro-combinatorial nature but its motivation
comes from geometry. The algebro-combinatorial properties studied
here should be useful in Schubert calculus associated with
orthogonal groups and the related degeneracy loci.

The computations of this paper are closely related to the ones in [LLT1];
we plan to develop this link in some future publication.

The algebro-combinatorial techniques, used in the present paper, 
are chosen to be
as elementary as possible. This should help the reader, with more
geometric and less algebro-combinatorial background, to read the paper. 
We mention however, that several results,
used in the proof of Theorem 9, in the appendix, are particular instances
of more general properties of Hall--Littlewood polynomials (see [LLT2]
and [LP2]).

Let us remark that there is also an interesting algebra and combinatorics of
``isobaric divided differences", with associated 
{\it Grothendieck polynomials} (cf. [FL]). 

\smallskip

This work has extensively used ACE ([V]) for explicit computations.

\smallskip

It is our pleasure and honor to dedicate
the present article to the mathematician whose recent work has illuminated
important connections between geometry and combinatorics. 

\bigskip\smallskip

\noindent
{NOTATION AND CONVENTIONS :}
\smallskip

A {\it vector} (of {\it length} $m$) is a sequence 
$[v_1,\ldots,v_m] \in \Bbb Z^m$. \hfill\break  
We will compare vectors of the same lengths, writing
$$
[v_1,\ldots,v_m] \subseteq [u_1,\ldots,u_m]
$$  
if $v_i \le u_i$  for all $i=1,\ldots,m$.

Given a vector $\alpha=[\alpha_1,\ldots,\alpha_m]$, we will write
$|\alpha|$ for the sum of its components.

A {\it partition} is an equivalence class of sequences 
$[i_1 \ge \cdots \ge i_m] \in \N^m$, where we identify the sequences 
$[i_1,\ldots, i_m]$ with $[i_1,\ldots,i_m,0]$. We denote the corresponding
partition by $I=(i_1,\ldots,i_m)$, by taking any representative
sequence.

A {\it part} of a partition $I$, is a nonzero component of any sequence
that represents $I$.  

The {\it length} of a partition $I$ is the number of its nonzero parts,
denoted $\ell(I)$. 

We call a partition {\it strict} if all its parts 
are different. 

We write $I\subseteq J$ for two partitions $I$ and $J$ 
(of possibly different lengths) if the same relation holds for any 
pair of the same length representing them. 

All operators act, in this paper, on their {\it left}.

Polynomials are usually treated as operators acting by multiplication.

\bigskip\medskip

\head
{\bf 1. Divided differences}
\endhead

Let $n$ be a fixed (throughout the paper)  positive integer.

The symmetric group (i.e. the Weyl group of type A) $\frak S_n$ 
is the group with generators $s_1,\ldots,\, s_{n-1}$
subject to the relations
$$
s_i^2=1, \quad \quad s_{i-1}\, s_i\, s_{i-1} =s_i s_{i-1}\, s_i
\quad ,\quad s_is_j =s_js_i \quad \forall i,j: |i-j|>1 \ .  \eqno (1.1) 
$$
We shall call $s_1,\ldots,\, s_{n-1}$ \ 
{\it simple transpositions} of $\S_n$. 

The hyperoctahedral group (i.e. the Weyl group of type B) $\frak B_n$ 
is an extension of $\S_n$ by an element $s_0$ such that
$$
s_0^2=1, \quad s_0\, s_1\, s_0\, s_1=s_1\, s_0\, s_1\, s_0,
\quad s_0\, s_i=s_i\, s_0 \quad \hbox{for} \quad i\ge 2.  \eqno (1.2) 
$$
The Weyl group $\frak D_n$ of type $D$ is the extension of  $\S_n$ by an element 
$\ss$ such that 
$$ \ss^2=1, \quad s_1 \ss\ =\ss\, s_1, \quad
\ss\, s_2\, \ss \, = s_2 \ss\, s_2, \quad   \ss\, s_i=s_i\, \ss 
\quad \hbox{for} \quad i> 2.   \eqno (1.3) 
$$
The group $\frak D_n$ can be thought as a subgroup of $\frak B_n$ by sending 
$\ss$ to  $s_0\, s_1\, s_0$. 

\medskip

The above three groups act on vectors of length  $n$ by
$$ [v_1,\ldots, v_n]\, s_i := [v_1,\ldots, v_{i-1},v_{i+1},v_i,
v_{i+2},\ldots,v_n] 
$$
$$ [v_1,\ldots, v_n]\, s_0 := [-v_1, v_2,\ldots,v_n] $$
$$ [v_1,\ldots, v_n]\, \ss := [-v_2,-v_1,\, v_3,\ldots,v_n]$$
The orbit of the vector $v=[1,\ldots,\, n]$ under $W=\S_n,\, \frak B_n$ or $\frak D_n$,
is in bijection with the elements of $W$, and we shall code 
each $w\in W$ by the vector $[1,\ldots,n]\, w$,
writing $\ov \imath$ instead of $-i$. 

\medskip

The three groups $W= \S_n,\, \frak B_n,\, \frak D_n$ also act on the ring of polynomials
in $n$ indeterminates $\Cal X=\{ x_1,\ldots, x_n\} $:
the simple transposition
$s_i$, $i\geq 1$, exchanges $x_i$ and $x_{i+1}$, $s_0$ sends $x_1$ to
$-x_1$,
$\ss$ sends $x_1$ to $-x_2$ and $x_2$ to $-x_1$, the action
being trivial in the non-listed cases. 
We shall denote by $f^w$ the image of a polynomial $f\in \Z[\Cal X]$  
under $w\in W$, and write $x^\alpha$, 
with $\alpha=[\alpha_1,\ldots,\alpha_n]\in \N^n$, 
for the monomial $x_1^{\alpha_1}\cdots x_n^{\alpha_n}$. 

\medskip

Let $\Pol$ be the ring of polynomials in $\X$
with coefficients in $\Z[{1\over 2}]$ (we need only division by 2).
For any $m\leq n$, let $\Sym(m|n-m)$ denote the subring of $\Pol$ 
consisting of polynomials invariant
under all $s_i$, $1\leq i\leq n-1$, $i\neq m$, and let 
$\Sym(n)=\Sym(n|0)=\Sym(0|n)$ be the ring of symmetric polynomials.
It contains as subrings $\Sym^B(n)$, the ring of polynomials invariant 
under $\frak B_n$, and   $\Sym^D(n)$, the invariants of $\frak D_n$.
It is easy to see that $\Pol$ is a free module over these different rings,
generated by $x^\alpha$, $\alpha\subseteq [n-1,\ldots,0]$ 
(or $\alpha\subseteq [0,\ldots,n-1]$)
over $\Sym(n)$,
by $x^\alpha$, $\alpha\subseteq [2n-1,2n-3,\ldots,1]$ 
(or $\alpha\subseteq [1, \ldots, 2n-3,2n-1]$)
over $\Sym^B(n)$, and
by $x^\alpha$, $\alpha\subseteq [2n-2,2n-4,\ldots,2,0]$ 
(or $\alpha\subseteq [0,2,\ldots,  2n-2]$) over $\Sym^D(n)$.

\medskip

The respective elements of maximal length in each of the groups are 
$$\omega := [n,\ldots,\, 1] \ \hbox{for}\ \S_n \ ,$$
$$ w_0^B := [\ov 1,\ldots,\, \ov n] \ \hbox{for}\ \frak B_n\ ,$$

$$ w_0^D := \left\{ \matrix [\ov 1,&\ldots,&\ov n] & n & \hbox{even} \cr
                [1,&\ov 2,&\ldots,&\ov n]&n &\hbox{odd} \cr
\endmatrix 
\right.
$$
for $\frak D_n$.
We shall also need the following element of $\frak D_n$:
$$ \vv := \omega w_0^D = 
\left\{ \matrix [\ov n,&\ldots,&\ov 1] & n & \hbox{even} \cr 
[n,&\ov{n-1},&\ldots,&\ov1] & n & \hbox {odd\, .} \cr
\endmatrix
\right.
$$

\smallskip

Relations between reduced decompositions in $W$ can be represented planarly.
By definition, a planar display will be identified with its reading from
left to right and top to bottom ({\it row--reading}). We shall also use
{\it column--reading}, that is, reading successive columns downwards,
from left to right.

For example, we will write 
$$\matrix 2 \cr 1 & 2 \cr \endmatrix  \equiv\  
\matrix 1 & 2 \cr & 1 \cr \endmatrix 
$$
for the following equality for simple transpositions:
$$ 
s_2\, s_1\, s_1 = s_1\, s_2\, s_1 \ .
$$

Suppose that a rectangle is filled 
row-wise from left to right, and column-wise from bottom to top
with consecutive numbers from $\{1,\ldots, n-1\}$. 

Then one easily checks that 
its row--reading and column--reading produce
two words which, interpreted as words in the $s_i$, are congruent
modulo the Coxeter relations.

Here is an example of such a congruence :
$$
\matrix 3456 \cr &2345 \cr &&1234 \cr \endmatrix \equiv
\matrix 3 \cr 2 \cr 1 \cr \endmatrix \cdot 
\matrix 4 \cr 3 \cr 2 \endmatrix \cdot \matrix 5 \cr 4 \cr3 \endmatrix
\cdot \matrix 6 \cr 5 \cr 4 \cr \endmatrix  
$$
the congruence class being conveniently denoted by the rectangle
$ \smatrix{3 & 4 & 5 & 6 \cr 2 & 3 & 4 & 5 \cr 1 & 2 & 3 & 4 \cr}$.

More generally, the planar arrays that we shall write, will have the
property that their row--reading and column--reading are congruent
modulo Coxeter relations \ (cf. [LS2,LS3] or [EG] for a ``jeu de taquin" 
on reduced decompositions). 
In this notation, one has, for any integers $a,b,c,d,k$ :
$1\leq a<b,c<d\leq n$, $a+d=b+c$, $k<d-b$, the congruence

$$ \matrix b+1 & \cdots & b+k \cr
            b  & \cdots & \cdots & \cdots & d \cr 
	               \vdots &     &        &      & \vdots \cr
		              a   & \cdots & \cdots & \cdots & c \cr  
\endmatrix				
\quad \equiv \quad   
\matrix \cr   b  & \cdots & \cdots & \cdots & d \cr
		\vdots &     &        &      & \vdots \cr
   			a   & \cdots & \cdots & \cdots & c \cr	
						& & c-k & \cdots & c-1 \cr
\endmatrix				
\eqno(1.4)
$$

\smallskip

It is convenient to work in the group algebra of $W=\S_n,\, 
\frak B_n$ or $\frak D_n$.
The works of Young and Weyl have stressed the r\^ole of the alternating sum
of elements of these groups.  Let, for $W=\S_n,\, \frak B_n$ or $\frak D_n$,
$$  
\Omega^W:= \sum_{w\in W} \, (-1)^{\ell(w)} w \, .   \eqno(1.5)  
$$
Using that $\frak B_n$ (resp. $\frak D_n$) is isomorphic to the semi-direct product
$\S_n \semi \Z_2^n$  (resp. $\S_n \semi \Z_2^{n-1}$), one obtains
the following factorizations in the group algebra:
$$ 
\Omega^{\frak B_n} = \Omega^{\S_n}\, \prod_{1\leq i \leq n} (1-\t_i)
= \prod_{1\leq i \leq n} (1-\t_i) \,\, \Omega^{\S_n} \, , 
\eqno(1.6) 
$$
$$ 
\Omega^{\frak D_n} = {1\over 2}\, \Omega^{\S_n}\, \Bigl(
\prod_{1\leq i \leq n} (1+\t_i) + \prod_{1\leq i \leq n} (1-\t_i)\Bigr)
$$    
$$
={1\over 2} \Bigl(
\prod_{1\leq i \leq n} (1+\t_i) +\prod_{1\leq i \leq n} (1-\t_i)\Bigr)\,\, 
\Omega^{\S_n}\, ,  
\eqno(1.7) 
$$
where $\t_1:=s_0$, and $\t_i= s_{i-1}\, \t_{i-1}\, s_{i-1}$ for $i>1$.
The elements $\Omega^W$, as operators on the ring of polynomials $\Pol$,
can be obtained from the cases of $\S_2$, $\frak B_1$, $\frak D_2$. 
To see this, we first need to define {\it simple divided differences}~:
$$
\Pol \ni f \mapsto f\, \d_i :=  (f-f^{s_i})/(x_i-x_{i+1}), 
 \quad i\ge 1 \,,
\eqno(1.8) 
$$ 
$$ 
\Pol \ni f \mapsto f\, \d_0 :=  (f-f^{s_0})/{(-x_1)}\ ,   
\eqno(1.9)
$$
$$ 
\Pol \ni f \mapsto f\, \dd := (f-f^{\ss} )/{(-x_1 - x_2)} \ .
\eqno(1.10)
$$

The $\d_i,\, \d_0,\, \dd$ satisfy the Coxeter relations (1.1)--(1.3), 
together with 
the relations
$$
\dd^2= 0 =   \d_i^2 \ \ \ \hbox{for} \ \ \ 0\leq i <n \ . 
\eqno(1.11)
$$
Therefore, to any element $w$ of the group $W$, there corresponds a 
{\it divided difference} $\d_w$. Any reduced decomposition
 \ $s_{i_1}\, s_{i_2}\, \cdots s_{i_{\ell}}= w$ \ of $w$ gives rise 
to a factorization \ $\d_{i_1}\, \d_{i_2}\, \cdots \d_{i_{\ell}}$ \ of 
$\d_w$ (cf. [BGG] and [D1,D2]).

We shall display divided differences planarly according to the same
conventions as for products of $s_i$'s.

For example, the divided difference 
$$
\d_0\d_1\d_2\d_3 \d_0\d_1\d_2\d_0\d_1
$$
will be displayed as
$$
\matrix \d_0 & \d_1 & \d_2 & \d_3 \cr
                             & \d_0 & \d_1 & \d_2 \cr
                                &  & \d_0 & \d_1 \cr
                              			\endmatrix
$$
As said before, the displays that we write have the property that their
row--reading is congruent to their column--reading, and thus the
preceding one encodes the equality
$$
\d_0 \d_1 \d_0 \d_2 \d_1 \d_0 \d_3 \d_2 \d_1 
= \d_0\d_1\d_2\d_3 \d_0\d_1\d_2\d_0\d_1\,.
$$

\medskip

We shall especially need the {\it maximal} divided differences
$\d_\omega$, $\d_{w_0^B}$, and $\d_{w_0^D}$. To describe them using alternating
sums of group elements, we define 
$$ 
\Delta:= \prod_{1\leq i<j\leq n} (x_i -x_j) =
	    x_1^{n-1}\cdots x_n^0\, \Omega^{\S_n} \ ,
\eqno(1.12)
$$
$$ 
\Delta^{B} := \prod_{i=1}^nx_i \prod_{n\ge i > j\ge 1} (x_i^2 -x_j^2)
	        =  {1\over 2^{n-1}} x^{[1,3,\ldots,2n-1]} \, \Omega^{\frak B_n} \ ,
\eqno(1.13)
$$
$$ \hbox{and} \ \ 
\Delta^{D} := \prod_{n\ge i>j \ge 1} (x_i^2 -x_j^2) =
	        {1\over 2^{n-1}}x^{[2,4,\ldots,2n-2]} \, \Omega^{\frak D_n} \ .
\eqno(1.14)
$$ 
 
The Weyl character formula for types A, B, and D can be written as 

\proclaim{Lemma 1} \ 
For each of the groups $W=\S_n,\, \frak B_n$ or $\frak D_n$, the alternating sum 
$\Omega^W$, as an operator on the ring of polynomials $\Pol$, 
is related to the maximal divided difference by
$$ 
\Omega^{\S_n} {1\over \Delta}  = \d_\omega\,,  \ \ \ \ \
\Omega^{\frak B_n}  {1\over \Delta^B} = (-1)^{n\choose 2}\, \d_{w_0^B}\,,  \ \ \
\hbox{and} \ \ \
\Omega^{\frak D_n}  {1\over \Delta^D} = (-1)^{n\choose 2} \, \d_{w_0^D} .
$$ 
\endproclaim
                                      
Indeed, all the operators in Lemma 1 commute with multiplication 
by polynomials 
which are invariant under $W$. Moreover, they decrease degree by the length
of the maximal element of the group.  Since $\Pol$ is a module over 
$\Sym^W(n)$ with a basis of monomials of degree strictly less than this
length, except for a single monomial, it remains only to check that
the actions of $\Omega$'s and $\partial$'s agree on this monomial, which
offers no difficulty. 

\bigskip\medskip

\head
{\bf 2. Bases of polynomial rings}   
\endhead 

The monomials mentioned in the previous section are not an appropriate basis, 
when interpreted in terms 
of cohomology classes for the flag variety.
Define, for the rest of this paper, the vector
$$
\rho:=[n-1,\ldots, 1,0]\, .
\eqno(2.1)
$$ 
Motivated by geometry, one defines recursively {\it Schubert polynomials} 
$Y_\alpha$, for any sequence $\alpha\in \N^n$, with $\alpha\subseteq \rho$, 
by
$$ 
Y_\alpha\, \d_i = Y_{\beta} 
 \ , \ \ \ \hbox{if} \ \ \alpha_i >\alpha_{i+1} \ \ ,  \eqno (2.2) 
$$
where 
$$
\beta = [\alpha_1,\ldots,\alpha_{i-1}, \alpha_{i+1},
\alpha_i-1, \alpha_{i+2},\ldots, \alpha_n]\,,
$$
starting from $Y_\rho= x^\rho$   (cf. [LS1], [M2]). 

In particular, if $\alpha\in \N^n$ is weakly decreasing, 
then $Y_\alpha$ is equal to the {\it monomial} \ $x^\alpha$\,. 

If, on the contrary, 
$\alpha_1\leq \cdots \leq \alpha_k $ and $\alpha_{k+1}=\cdots =\alpha_n=0$, 
for some $k\le n$,
then 
$Y_\alpha$ coincides with the {\it Schur polynomial} $s_\l(x_1,\ldots,x_k)$,  
where $\l = (\alpha_k,\ldots, \alpha_1)$.

\medskip

\noindent
CONVENTION: \ Let $\alpha \in \Bbb N^k$.
Then we shall write $Y_{\alpha}$ for $Y_{\alpha,0,\ldots,0}$.

\smallskip

We also record, for later use, the following equality: for  
$\alpha = [\alpha_1, \ldots ,\alpha_k] \in \Bbb N^k$, 
$$
Y_{\alpha} \, x_1 \cdots x_k = Y_{[\alpha_1+1,\ldots, \alpha_k+1]}.
\eqno(2.3)
$$  

\smallskip

On $\Pol$ there is a scalar product:
$$
( \ , \ ): \Pol \times \Pol \to \Sym(n),
$$
defined for $f, g\in \Pol$ by
$$
(f,g): = f g \, \d_{\omega}. \eqno(2.4)
$$

There exists an involution
\footnote {This involution is: \ $\hbox {code}(w)=
\alpha \mapsto \alpha'=\hbox{code} (w \omega)$ \ (cf. [M2]).}
\ $\alpha \mapsto \alpha'$ \  
such that
$$
\bigl(Y_{\alpha}^\omega , Y_{\beta'}\bigr) = (-1)^{|\alpha|} 
\delta_{\alpha \beta}\,. 
\eqno(2.5)
$$
Moreover, when $\alpha,\beta \subseteq \rho$ are such that $|\alpha| + 
|\beta| = |\rho|$, then one has
$$
\bigl(Y_{\alpha}, Y_{[n-1-\beta_1, n-2-\beta_2,\ldots,0-\beta_n]}\bigr)=
\delta_{\alpha \beta}.
\eqno(2.6)
$$

We also will need \ $\Q$-polynomials \ of [PR].
We set $\Q_i:= e_i=e_i(\X)$, the $i$-th elementary
symmetric polynomial in $\X$\,.
Given two nonnegative integers $i \ge j$, we adapt Schur's definition of
his $Q$-functions by putting
$$
\Q_{i,j}:=\Q_i \Q_j
+2\sum\limits^j_{p=1}(-1)^p\Q_{i+p} \Q_{j-p}\,.
\eqno(2.7)
$$
Given any partition $I=(i_1,\ldots,i_k)$, where we can assume $k$ to be even,
we set
$$
\Q_I:= \hbox{Pfaffian} (M)\,,
\eqno(2.8)
$$
where $M=(m_{p,q})$ is the $k\times k$ skew-symmetric matrix
with $m_{p,q}= \Q_{i_p,i_q}$  for $1 \le p < q \le k$.

\smallskip

Equivalently, for any partition $I=(i_1\ge i_2\ge\ldots\ge i_{\ell}>0)$\,,
the polynomial $\Q_I=\Q_I(X)$ is defined recurrently
on $\ell$ by putting for odd $\ell$\,,
$$
\Q_I:=\sum\limits^{\ell}_{j=1}(-1)^{j-1}
\Q_{i_j} \Q_{(i_1,\ldots,i_{j-1},i_{j+1},\ldots,i_{\ell})}
\eqno(2.9)
$$
and for even $\ell$\,,
$$
\Q_I:=\sum\limits^{\ell}_{j=2}(-1)^j
\Q_{i_1,i_j} \Q_{(i_2,\ldots,i_{j-1},i_{j+1},\ldots,
i_{\ell})}\,.
\eqno(2.10)
$$
For any positive integer $k$, let $\rho(k)$ denote the partition 
$$
\rho(k):=(k,k-1,\ldots, 1)\,.
\eqno(2.11)
$$
The ring $\Sym(n)$ is a free module over the ring of polynomials
symmetric in $x_1^2,\ldots,x_n^2$, with a basis 
provided by the $\Q_I(\X)$, where $I\subseteq \rho(n)$ ranges over 
strict partitions. 

As functions of $x_1,\ldots, x_m$, the $\Q$-polynomials 
can also be defined recursively
by induction on $m$, involving now all partitions without restriction,
as for Hall-Littlewood polynomials: for any strict partition $I$, one has
$$
\Q_I(x_1,\ldots,x_m) = \sum_{j=0}^{\ell(I)} x_m^j \bigr( \sum_{|I|-|J|=j}
\Q_J(x_1,\ldots,x_{m-1})\bigl),
\eqno(2.12)
$$
where the sum is over all (i.e. not necessarily strict) 
partitions $J\subseteq I$
such that $I/J$ has at most one box in every row (cf. [PR, Prop. 4.1]).
Moreover, given a
partition $I'=(\ldots,i,j,j,k,\ldots)$ and denoting $I=(\ldots,i,k,\ldots )$,
one has the {\it factorization property}
$$
\Q_{I'}=\Q_{j,j} \Q_I.
\eqno(2.13)
$$

We define, for a strict partition $I$,  
$$
\P_I:= 2^{-\ell(I)}\, \Q_I\,.
\eqno(2.14)
$$
The ring $\Sym(n)$ is a free module over $\Sym^D(n)$ with a basis 
provided by the $\P_I$, where $I$ ranges over strict partitions 
contained in $\rho(n-1)$.

Now we will need the following divided difference:
$$
\d_v 
=(\dd \d_2\cdots \d_{n-1}\, \d_1\cdots \d_{n-2}) 
\cdots (\dd \d_2\d_3\, \d_1\d_2)\, \dd  \ \ \ \ \hbox{$n$ even}   
\eqno(2.15)
$$
and
$$
\d_v
=(\dd \d_2\cdots \d_{n-1}\, \d_1\cdots \d_{n-2})
         \cdots (\dd \d_2\d_3\d_4\, \d_1\d_2\d_3)\, (\dd \d_2\d_1) 
	 \ \ \hbox{ $n$ odd}  
\eqno(2.16)
$$

Denote by
$$
\langle \ ,\ \rangle : \Sym(n) \times \Sym(n) \to Sym^D(n) 
$$
the scalar product defined for $f, g\in Sym(n)$ by
$$
\langle f,g \rangle := f g \, \d_v  .
\eqno(2.17)
$$
For strict partitions $I,J\subseteq \rho(n-1)$, one has
$$
\langle \P_I,\P_{\rho(n-1) \smallsetminus J} \rangle 
= (-1)^{n\choose 2}\, \delta_{I J}, 
\eqno(2.18)
$$
where $\rho(n-1) \smallsetminus I$ is the strict partition whose parts 
complement the parts of $I$ in $\{n-1,n-2,\ldots, 1\}$ \ (cf. [PR]).

Consequently, the polynomial ring $\Pol= \Bbb Z [{1 \over 2}][x_1,\ldots, x_n]$
is a free $\Sym^D(n)$-module with a basis $Y_{\alpha}\, \P_I$, where
$\alpha$ ranges over subsequences contained in $\rho$ and $I$
runs over all strict partitions contained in $\rho(n-1)$. 
Note that the element of maximal degree of this basis is
$x^{\rho} \P_{\rho(n-1)}$. Let
$$
[ \ , \ ]: \Pol \times \Pol \to \Sym^D(n)
$$
be a scalar product, defined for $f,g \in \Pol$ by
$$
[f,g]:= f g \,\d_{w_0^D}.
\eqno(2.19)
$$ 
One has, for $\alpha, \beta \subseteq \rho$ and strict partitions $I,J\subset
\rho(n-1)$,
$$
\bigl[Y_{\alpha}^{\omega}\, \P_I , 
Y_{\beta'}\, \P_{\rho(n-1)\smallsetminus J} \bigr] = 
(-1)^{|\alpha|+{n \choose 2}}\, \delta_{\alpha \beta}\, \delta_{I J}\,.
\eqno(2.20)
$$ 
(See (2.5).)

\medskip

Let $\Y=\{ y_1,\ldots,y_n\}$ be a second set of indeterminates of
cardinality $n$. 
The symbol $\equiv$ will mean: \ ``congruent modulo the ideal
generated by the relations \hfill\break
$ f(x_1^2,\ldots, x_n^2) =
f(y_1^2,\ldots, y_n^2)$, $f\in Sym(n)$,   
together with  $x_1\cdots x_n = y_1\cdots y_n$".

\medskip

Following Fulton [F2,F3], define 
$$
\aligned
F(\X,\Y) &:=\ |\P_{n+j-2i}(\X)+\P_{n+j-2i}(\Y)|_{1\le i,j\le n-1} \\
& \\
&=\ {\left|\ \CD
\P_{n-1}(\X)+\P_{n-1}(\Y) & 0 &.\;.\;. & \\
\P_{n-3}(\X)+\P_{n-3}(\Y) &\ \P_{n-2}(\X)+\P_{n-2}(\Y) &\ .\;.\;.\ & \\
\vbox{\offinterlineskip\hbox{.}\vskip3pt\hbox{.}\vskip3pt\hbox{.}}
&\vbox{\offinterlineskip\hbox{.}\vskip3pt\hbox{.}\vskip3pt\hbox{.}}
&\vbox{\offinterlineskip\hbox{.}\vskip3pt\hbox{\ \ .}\vskip3pt
\hbox{\ \ \ \ .}} & \\
&& 1 &&  \P_1(\X)+\P_1(\Y) \\
\endCD\ \right|}\,. \\
\endaligned
\eqno(2.21)
$$

Following [PR], define
$$ 
\P(\X, \Y):= \sum \P_I(\X)\, \P_{\rho(n-1) \smallsetminus I}(\Y) \ ,  
\eqno(2.22)
$$
where the summation is over all strict partitions $I\subseteq \rho(n-1)$.
The reasoning in [LP1, Sect.2] made for case $C_n$ adapts to case $D_n$ 
and furnishes:

\proclaim {Proposition 2} \ We have 

\noindent
(i) \ $$F(\X,\Y) \equiv \P(\X,\Y)\,. \eqno(2.23)$$ 

\noindent
(ii) \ For every $w\in \frak D_n\smallsetminus \S_n$,
$$\P(\X^w, \X)=0\,,
\eqno(2.24)
$$
and for every $w\in \S_n$,
$$
\P(\X^w, \X)= \P(\X, \X)= s_{\rho(n-1)}(\X)\,.
\eqno(2.25)
$$

\noindent
(iii) \ For every $f\in \Sym(n)$, 
$$
\langle f(\X)\,, F(\X,\Y) \rangle \equiv (-1)^{n\choose 2}\, f(\Y) \ . 
\eqno(2.26)
$$

\noindent
(iv) \ For every $f\in \Pol$,
$$
\bigl[ f(\X)\,, \prod_{n \ge i>j \ge 1}(x_i-y_j)\, F(\X,\Y) \,\bigr] 
\equiv f(\Y)\,.
\eqno(2.27)
$$
\endproclaim
In other words, $F(\X,\Y)$ is a reproducing kernel for the scalar product
$\langle \ , \ \rangle$, and $\prod_{i>j}(x_i-y_j)\, F(\X,\Y)$ is a reproducing
kernel for $[ \ , \ ]$. One can show that the ``vanishing property"
{\it (ii)} characterizes $\P(\X,\Y)$ up to $\equiv$. The congruence
{\it (i)} can be also derived from geometry by comparing the classes
of diagonals in flag bundles associated with $SO(2n)$ given in
[F2,F3] and [PR] (see also [G]).

\bigskip\medskip

\head 
{\bf 3. Vertex Operators}
\endhead

In this section we shall mainly make computations using
the following two divided differences:

\proclaim {Definition 3} \ For $k\le n$, we set
$$
\nabla_k^B(n):= (\d_0 \d_1 \cdots \d_{n-1}) \cdots (\d_0 \d_1 
\cdots \d_{n-k}) \,.
\eqno(3.1)
$$
For $k\le n/2$, we put
$$
\nabla_k^D(n) :=(\dd \d_2\cdots \d_{n-1} \d_1 \d_2 \cdots \d_{n-2})
             \cdots  (\dd \d_2\cdots \d_{n-2k+1} \d_1 \d_2 \cdots \d_{n-2k}).
\eqno(3.2)
$$
\endproclaim

We shall need the following fact from [LP1], quoted in the appendix:

\proclaim{Fact 4} \ Let $k\leq n$ and let 
$\alpha =[\alpha_1 \leq \cdots \leq \alpha_k] \in \N^{k}$ with 
$\alpha_k \le n-k$.
Suppose that $I\subseteq \rho(n)$ is a strict partition. 
Then the image of $\Q_I \, Y_\alpha$ under $\nabla_k^B(n)$ is 0 
unless $n-0-\alpha_1,\ldots, n-(k-1) -\alpha_k$ are
parts of $I$. In this case, the image is 
$(-1)^{k(n-1)+s} 2^k\Q_J$, 
where $J$ is the strict partition with parts
$$\{i_1,\ldots, i_{\ell(I)} \} \smallsetminus
\{ n-0-\alpha_1,\ldots, n-(k-1)-\alpha_k \},
$$
and $s$ is the sum of positions of the parts erased in $I$.
\endproclaim 

\bigskip

\noindent
{\bf Example 5.} \ For $n=7$ and $k=2$\,, we have
$$
\aligned
&\Q_{(5,4,3,2,1)} \, Y_{[2,5]} \, \nabla_2^B(7) \cr
&=\Q_{(5,4,3,2,1)} \, Y_{[2,5]}\, (\d_0 \d_1 \d_2 \d_3 \d_4 \d_5 \d_6) 
(\d_0 \d_1 \d_2 \d_3 \d_4\d_5)= 4 \Q_{(4,3,2)}\,, \cr
\endaligned
$$
and for $k=3$\,, we have
$$
\Q_{(7,5,4,3,1)} \, Y_{[2,3,4]} \nabla_3^B(7) = - 8 \Q_{(7,4)}\,.
$$

\bigskip

The following result establishes a basic relation between the $\nabla^D$'s and
$\nabla^B$'s:

\proclaim{Proposition 6} \quad  Let $k$ be a positive integer.
As operators on $Sym(2k)$,
$$
\nabla_k^D(2k) = x_1\cdots x_{2k} \, \nabla_{2k}^B(2k)
   +x_1\cdots x_{2k-1} \, \nabla_{2k-1}^B(2k) \,.
\eqno(3.3)
$$
\endproclaim
Before proving (3.3), we illustrate it by the following examples:

\bigskip
\noindent
{\bf Example 7.} As operators on $\Sym(2)$,
$$ 
\nabla_1^D(2) = \dd = x_1x_2 \d_0\d_1\d_0 + x_1 \d_0\d_1.
$$
As operators on $\Sym(4)$,
$$
\aligned 
\nabla_2^D(4) &=  (\dd\d_2\d_3\d_1\d_2)\dd= \cr 
   &x_1x_2x_3x_4\, (\d_0\d_1\d_2\d_3)(\d_0\d_1\d_2)(\d_0\d_1)\d_0
    +x_1x_2x_3\, (\d_0\d_1\d_2\d_3)(\d_0\d_1\d_2)(\d_0\d_1). 
\endaligned
$$

\medskip
\noindent
The RHS of the last equation is depicted planarly as
$$
x_1 x_2 x_3 x_4 \matrix \d_0 & \d_1 & \d_2 & \d_3 \cr
                             & \d_0 & \d_1 & \d_2 \cr
                                &  & \d_0 & \d_1 \cr
                                     &  &  & \d_0 \cr 
			\endmatrix
+x_1 x_2 x_3 \matrix \d_0 & \d_1 & \d_2 & \d_3 \cr
                             & \d_0 & \d_1 & \d_2 \cr
                                &  & \d_0 & \d_1 \cr
                              			\endmatrix
$$

\bigskip

\demo {Proof of the proposition} \
In this proof, let $\X:=\{x_1,\ldots,x_{2k} \}$.
Both sides of (3.3) are $Sym^B(2k)$-linear.
The operator $\nabla_k^D(2k)$ sends all $\Q_I(\X)$, 
$I\subseteq \rho(2k-1)$ to 0,
except for $\Q_{\rho(2k-1)}$ which is sent to $(-1)^k 2^{2k-1}$ (cf.(2.18)). 
Thus $\nabla_k^D(2k)$ annihilates all $\Q_I(\X)$, $I\subseteq \rho(2k)$, 
except for $I=\rho(2k-1)$ which is sent to $(-1)^k 2^{2k-1}$ and 
$I=\rho(2k)$ which is sent to $(-1)^k 2^{2k-1} x_1 \cdots x_{2k}$.

The action of 
$$
x_1\cdots x_{2k} \nabla_{2k}^B(2k)
$$ 
is given by Fact 4.  Only $\Q_{\rho(2k)}(\X)$ survives and is sent to 
$(-1)^k 2^{2k} x_1\cdots x_{2k}$. 

We will now calculate the action of 
$$
x_1\cdots x_{2k-1} \, \nabla_{2k-1}^B(2k)
$$ 
on the $\Q_I(\X)$, where 
$I\subseteq \rho(2k)$ is a strict partition.
We set, temporarily in this proof, 
$$
\nabla:=\nabla_{2k-1}^B(2k) \ \ \ \ \hbox{and} \ \ \ \  
\nabla':=\nabla_{2k-1}^B(2k-1)\,,
$$ 
so that 
$$
\nabla=\nabla' \d_{2k-1} \cdots \d_1.
$$
Let $\X':=\{x_1,\ldots,x_{2k-1}\}$. 
We decompose $\Q_I(\X)$ as a sum of products of powers of $x_{2k}$
times some $\Q_J(\X')$, according to the formula (2.12):
$$
\Q_I(\X)= \sum \Q_J(\X')\, x_{2k}^{m_J} \ .
$$
Let $\overline J$ be the strict partition obtained from  a partition 
$J$ by subtracting all the pairs of equal parts. 
We have three cases to examine:
\smallskip

\noindent
1. Let $i_1\le 2k-2$. Then for each $J$, $|\overline J| + m_J + 2k-1 < 
\hbox {deg} \nabla$, and hence 
$$
x_1\cdots x_{2k-1} \Q_I(\X) \nabla = 0.
$$
 
\smallskip

\noindent
2. Let $i_1=2k-1$. For degree reasons, $\Q_I(\X)x_1\cdots x_{2k-1}
\nabla \ne 0$ is possible
only if $I=\rho(2k-1)$ ($I$ being a strict partition).

\noindent
\proclaim {Claim} \ We have
$$
x_1\cdots x_{2k-1} \Q_J(\X') \nabla' \ne 0
\eqno(3.4)
$$
only if $J = \overline J = \rho(2k-2)$.
\endproclaim

Indeed, suppose first that $j_1=2k-1$. Then
$$
x_1\cdots x_{2k-1} \Q_J(\X') = \Cal P \cdot \Q_H(\X'),
$$
where $\Cal P$ is a polynomial symmetric in $x_1^2,\ldots,x_{2k-1}^2$,
and a strict partition $H$ has no part equal to $2k-1$. Since this
expression is annihilated by $\nabla'$, we cannot have (3.4).
So, for degree reasons, (3.4) holds only if $j_1=2k-2$.
Suppose now that $j_2=2k-2$.
We get
$$
x_1\cdots x_{2k-1} \Q_J(\X') = \Cal P \cdot \Q_H(\X'),
$$
where $\Cal P$ is a polynomial symmetric in $x_1^2,\ldots,x_{2k-1}^2$,
and a strict partition $H$ has no part equal to $2k-2$. Since this expression
is annihilated by $\nabla'$, we cannot have (3.4). So, for degree
reason, (3.4) holds only if $j_2=2k-3$. Continuing this way,
we get the claim.

\smallskip

For $J=\rho(2k-2)$, we compute
$$
x_1\cdots x_{2k-1} \Q_J(\X') x_{2k}^{2k-1} \nabla
=\Q_{\rho(2k-1)}(\Cal X')\nabla' x_{2k}^{2k-1}\, \d_{2k-1} \cdots \d_1 =
(-1)^{k-1} 2^{2k-1}.
$$
\smallskip

\noindent
3. Let $i_1=2k$. Then $(i_2,\ldots ) \subseteq \rho(2k-1)$. We have
$$
\aligned
& x_1\cdots x_{2k-1} \Q_I(\X) \nabla = x_1^2 \cdots x_{2k-1}^2 x_{2k}
\Q_{(i_2,\ldots)}(\X) \nabla \cr
=& \Q_{(i_2,\ldots)}(\X) \nabla' x_1^2 \cdots x_{2k-1}^2 x_{2k}\, 
\d_{2k-1} \cdots \d_1 \ .
\endaligned
$$
Now for $H\subseteq (i_2,\ldots) \subseteq \rho(2k-1)$, $\Q_H(\X')\nabla' \ne 0$
iff $H=\rho(2k-1)$, so iff $(i_2,\ldots)=\rho(2k-1)$. We have
$$
\Q_{\rho(2k-1)}(\Cal X) \nabla' =\Q_{\rho(2k-1)}(\Cal X')\nabla' 
= (-1)^k 2^{2k-1}
$$
and
$$
x_1^2 \cdots x_{2k-1}^2 x_{2k}\, \d_{2k-1} \cdots \d_1 = - x_1 \cdots x_{2k}.
$$
Summarizing,
$$
\Q_I(\Cal X) x_1 \cdots x_{2k-1} \nabla \ne 0
$$ 
only if $I=\rho(2k-1)$, when we get $(-1)^{k-1} 2^{2k-1}$; 
or $I=\rho(2k)$, when
we get $(-1)^{k-1} 2^{2k-1} x_1\cdots x_{2k}$.

Finally, comparing the computed values of the $\Q_I(\X)$ under the operators:
$$
\nabla_k^D(2k)\,, \ \ \ \ \ x_1\cdots x_{2k}\,  \nabla_{2k}^B(2k), 
\ \ \hbox {and} \ \ \ \ \ x_1\cdots x_{2k-1} \, \nabla_{2k-1}^B(2k),
$$
that are possibly nonzero only for $I=\rho(2k)$ and $\rho(2k-1)$, 
we get the desired formula (3.3). (Note that we have also used the equality 
$2^{p-1}=2^p - 2^{p-1}$.)
\qed
\enddemo

\medskip

\proclaim{Corollary 8} \ Let $k$ be a positive integer such that $k \leq n/2$.
As operators on the ring \ $\Sym(2k\, |\, n-2k)$,   
$$
\nabla_k^D(n) = 
x_1\cdots x_{2k} \, \nabla_{2k}^B(n) 
+x_1\cdots x_{2k-1} \, \nabla_{2k-1}^B(n) \, \d_1 \cdots \d_{n-2k}\,.
\eqno(3.5)
$$
\endproclaim

This property is obtained from Proposition 6 by composing the  
expression for the operator $\nabla_k^D(2k)$ with the divided difference
$$
(\d_{2k} \cdots \d_{n-1}) \cdots (\d_2 \cdots \d_{n-2k+1})
(\d_1\cdots \d_{n-2k}) = \matrix 
\d_{2k} & \cdots & \d_{n-1} \cr
\vdots & & \vdots \cr
\d_2 & \cdots & \d_{n-2k+1} \cr
\d_1 & \cdots & \d_{n-2k} \cr
\endmatrix
$$ 

\bigskip

In the proof of Theorem 11, 
we will need the following supplement to Fact 4:

\proclaim{Theorem 9}  
\ Let $k\leq n$ and let 
$\alpha =[\alpha_1 \leq \cdots \leq \alpha_k] \in \N^{k}$ with 
$\alpha_k=n-k+1$.
Suppose that $I\subseteq \rho(n)$ is a strict partition. 
Then the image of \ $\Q_I \, Y_\alpha$ \  
under $\nabla_k^B(n)$ is 0 
unless $\ell(I)\not\equiv n$ (mod $2$) and 
$n-0-\alpha_1,\ldots, n-(k-2) -\alpha_{k-1}$ are
parts of $I$. 
In this case, the image is 
$(-1)^{(k-1)(n-1)+1+s} 2^k\Q_J$, 
where $J$ is the strict partition with parts
$$\{i_1,\ldots, i_{\ell(I)} \} \smallsetminus
\{ n-0-\alpha_1,\ldots, n-(k-2)-\alpha_{k-1} \} \,,
$$
and $s$ is the sum of positions of the parts erased in $I$.
\endproclaim 

The proof of this theorem will be given in the appendix.

\bigskip
\noindent
{\bf Example 10.}  

(i) For $n=5$ and $k=1$, we have 
$$
x_1^5 \Q_{(5,3,2,1)} \d_0 \d_1 \d_2 \d_3 \d_4 = - 2 \Q_{(5,3,2,1)}
\ \ \ \hbox{and} \ \ \  x_1^5 \Q_{(5,2,1)} \d_0 \d_1 \d_2 \d_3 \d_4 = 0.
$$

\noindent
(ii) For $n=7$ and $k=2$, we have
$$
\Q_{(7,6,4,1)}\, Y_{[1,6]}\, \nabla_2^B(7) 
= - 4 \Q_{(7,4,1)}
$$
and
$$
\Q_{(7,6,4,3,1)}\, Y_{[1,6]}\, \nabla_2^B(7) = 0.
$$

\noindent
(iii) For $n=7$ and $k=4$, we have 
$$
\Q_{(7,6,4,3,2,1)}\, Y_{[1,2,2,4]}\, \nabla_{4}^B(7) 
=16 \Q_{(7,2,1)}
$$
and
$$
\Q_{(7,6,4,3,2,1)}\, Y_{[1,3,4,4]}\, \nabla_{4}^B(7) 
=-16 \Q_{(7,4,2)}.
$$

\bigskip

The following theorem is the main result of the present paper:
 
\proclaim{Theorem 11} \quad Let $k$ be a positive integer such that  
$k \leq n/2$.
Suppose that $I\subseteq \rho(n-1)$ is a strict partition. Let 
$\alpha = [\alpha_1 \leq \alpha_2\leq \cdots \leq \alpha_{2k}] 
\in \N^{2k}$ with $\alpha_{2k} \leq n-2k$.
Then the image of $\P_I \, Y_{\alpha}$ under $\nabla_k^D(n)$
is $0$ unless all the integers
$n-1 -\alpha_1,\ldots,\, n-2k -\alpha_{2k}$ belong to 
$\{ i_1,\ldots,i_{\ell(I)},0 \}$. 
In this case, the image is $(-1)^s \P_J$, where 
$J$ is the strict partition with parts 
$$
\{ i_1,\ldots,i_{\ell(I)} \} \smallsetminus 
\{ n-1-\alpha_1,\ldots, n-2k -\alpha_{2k} \} \, .
$$
Moreover, let $s'$ be the sum of positions of the parts erased in $I$,
and $s'':=\ell(I)+1$. Then $s=s'$  if  $\alpha_{2k} < n-2k$, and
$s=s' + s''$ if $\alpha_{2k} = n-2k$. \footnote{It is convenient
to treat here $0=n-2k-\alpha_{2k}$ as an ``extra part" 
of $I$, and take $s$ to be the sum of positions 
of all the parts erased in $I$, including the extra part.}
\endproclaim

\bigskip

\noindent
{\bf Example 12.} \ (i) For $n=7$ and $k=1$, we have
$$
\P_{(5,4,3,2,1,0)}\, Y_{[1,3]} \nabla_1^D(7)= 
\P_{(5,4,3,2,1,0)}\, Y_{[1,3]}\, \dd \d_2 \d_3 \d_4 \d_5 \d_6 
\d_1 \d_2 \d_3 \d_4 \d_5 = - \P_{(4,3,1)}
$$
and
$$
\P_{(6,4,3,2,1,0)}\, Y_{[2,5]} \, \nabla_1^D(7)= \P_{(6,3,2,1)}.
$$
\smallskip

\noindent
(ii) For $n=7$ and $k=2$, we have
$$
\P_{(6,5,4,3,2,1,0)} \, Y_{[1,1,1,2]}\, \nabla_2^D(7) = - \P_{(6,2)}
$$
and
$$
\P_{(6,5,4,3,2,1,0)} \, Y_{[1,1,1,3]}\, \nabla_2^D(7) = \P_{(6,2,1)}\,.
$$

\bigskip

\demo{Proof of the theorem} \ To compute the action of $\nabla_k^D(n)$, 
one uses its decomposition into a sum of two operators, given in Eq.(3.5). 

The image of \ $\Q_I(\X)\,  Y_\alpha$ \  
under the first operator 
$$
\Omega_1 := x_1\cdots x_{2k} \, \nabla_{2k}^B(n)
$$ 
is given by Fact 4 combined with Eq.(2.3) if \ $\alpha_{2k}<n-2k$\,, 
and by Theorem 9 combined with (2.3) in the case \ $\alpha_{2k}=n-2k$.

Since $x_{2k}$ appears in $ Y_\alpha$, the same results, however,  
{\it do not} directly furnish the value of \ $\Q_I(X_n)\,  Y_\alpha $ \  
under the second operator 
$$
\Omega_2 := x_1\cdots x_{2k-1} \, \nabla_{2k-1}^B(n)\,\d_1\cdots \d_{n-2k}\,.
$$ 
To end this computation, we proceed as follows.
For simplicity of indices, let us take temporarily \ $n=7$ and $k=2$. 
Suppose that $\alpha= [\alpha_1\le \alpha_2\le \alpha_3 \le \alpha_4] 
\in \N^4$ is such that $\alpha_4 \le 3$.
We want to compute
$$
 \Q_I \, Y_\alpha \, x_1x_2x_3\, 
  \matrix \d_0 &\d_1 &\d_2 &\d_3 &\d_4 &\d_5 &\d_6 \cr
     & \d_0 &\d_1 &\d_2 &\d_3 &\d_4 &\d_5 \cr
      && \d_0 &\d_1 &\d_2 &\d_3 &\d_4 \cr
       &&&& \d_1 &\d_2 &\d_3  \cr 
\endmatrix 
$$
Now, thanks to the relations (1.4), one has 

\smallskip

$$
\matrix  \d_0 &\d_1 &\d_2 &\d_3 &\d_4 &\d_5 &\d_6 \cr
          & \d_0 &\d_1 &\d_2 &\d_3 &\d_4 &\d_5 \cr
	  && \d_0 &\d_1 &\d_2 &\d_3 &\d_4 \cr 
         &&&& \d_1 &\d_2 &\d_3 \cr 
\endmatrix
$$

$$
=\left( \matrix \d_0 &\d_1 &\d_2\cr 
	  & \d_0 &\d_1\cr &&\d_0 \cr
\endmatrix	  
\right) \,
\left( \matrix \d_3 &\d_4 &\d_5 &\d_6 \cr
\d_2 &\d_3 &\d_4 &\d_5 \cr
\d_1 &\d_2 &\d_3 &\d_4\cr
&\d_1 &\d_2 &\d_3 \cr 
\endmatrix 
\right) 
= \left(\matrix \d_0 &\d_1 &\d_2\cr & \d_0 &\d_1\cr 
	  &&\d_0\cr
\endmatrix
\right)
\left(
\matrix \d_4 &\d_5 &\d_6 \cr 
	  \d_3 &\d_4 &\d_5 &\d_6 \cr\d_2 &\d_3 &\d_4 &\d_5 \cr
	  \d_1 &\d_2 &\d_3 &\d_4\cr
	  \endmatrix
	  \right)
$$

\medskip

$$	   
=\left(\matrix \d_0 &\d_1 &\d_2\cr & \d_0 &\d_1\cr &&\d_0\cr
	  \endmatrix
	  \right)  
	  (\d_4\d_5\d_6)
	  \left(\matrix \d_3 &\d_4 &\d_5 &\d_6 \cr
	  \d_2 &\d_3 &\d_4 &\d_5 \cr
	  \d_1 &\d_2 &\d_3 &\d_4 \cr
	  \endmatrix
	  \right)
$$
\smallskip

\noindent
Since $\d_4\d_5\d_6$ commutes with the divided differences on its left,
the last expression is rewritten as

\smallskip

$$ 
\d_4\d_5\d_6 \, 
\matrix \d_0 &\d_1 &\d_2 &\d_3 &\d_4 &\d_5 &\d_6 \cr
	     & \d_0 &\d_1 &\d_2 &\d_3 &\d_4 &\d_5 \cr 
	     && \d_0 &\d_1 &\d_2 &\d_3 &\d_4 \cr 
	     \endmatrix
= \d_4\d_5\d_6\, \nabla_3^B(7)
$$
\smallskip

\noindent

\noindent
Since $\d_4\d_5\d_6$ commutes with $x_1x_2x_3$ and $\Q_I= 
\Q_I(x_1,\ldots,x_7)$, the polynomial to be computed is equal to 

$$ 
Y_{\alpha} \, \d_4\d_5\d_6 \,   \Q_I \,  x_1x_2x_3\, \nabla_3^B(7)\,.  
$$
However, using (2.2), the image of $Y_{[\alpha_1, \alpha_2, \alpha_3, 
\alpha_4]}$ under $\d_4\d_5\d_6$ is 
$$
Y_{[\alpha_1, \alpha_2, \alpha_3, 0,0,0, \alpha_4-3]}
= Y_{[\alpha_1,\alpha_2, \alpha_3]}
$$ 
if $\alpha_4=3$, and 0 otherwise. Hence, by (2.3), the polynomial 
to be computed is equal to 
$$
Y_{[\alpha_1+1,\alpha_2+1,\alpha_3+1]} \Q_I \nabla_3^B(7)\,.
$$

In general, arguing along these lines, we evaluate
$$
\Q_I\, Y_{\alpha}\, x_1 \cdots x_{2k-1} \, \nabla^B_{2k-1}(n)
\, \d_1 \d_2 \cdots \d_{n-2k}\,,
$$
where $\alpha = [\alpha_1\le \alpha_2 \le \cdots \le \alpha_{2k}]\in \N^{2k}$
is such that $\alpha_{2k}\le n-2k$. By the relations (1.4), this amounts
to evaluating
$$
\Q_I\, Y_{\alpha}\, x_1 \cdots x_{2k-1} \,  \d_{2k} \d_{2k+1}\cdots 
\d_{n-1}\,\nabla^B_{2k-1}(n)\,.
$$
Since $\d_{2k} \d_{2k+1}\cdots \d_{n-1}$ 
commutes with $x_1 \cdots x_{2k-1}$
and $\Q_I = \Q_I(x_1,\ldots,x_n)$, the polynomial to be computed is equal to 
$$
Y_{\alpha}\, \d_{2k} \d_{2k+1}\cdots \d_{n-1} \Q_I x_1 \cdots x_{2k-1} \, 
\nabla^B_{2k-1}(n)\,.
$$
However, using (2.2), the image of $Y_{\alpha}$ under  
$\d_{2k} \d_{2k+1}\cdots \d_{n-1}$ is
$$
Y_{[\alpha_1,\ldots,\alpha_{2k-1},0^{n-2k},\alpha_{2k}-(n-2k)]}\,.
\eqno(3.6)
$$
The expression (3.6) equals $0$ unless $\alpha_{2k} =n-2k$,
when it is equal to
$Y_{[\alpha_1,\ldots,\alpha_{2k-1}]}$\,. Hence, by (2.3), the polynomial 
to be computed is equal to 
$$
\Q_I Y_{[\alpha_1+1,\ldots,\alpha_{2k-1}+1]} \, \nabla_{2k-1}^B(n)\,.
$$
Since $\alpha_{2k-1}+1 \le n-(2k-1)$, Fact 4 provides the end of the 
computation with the operator $\Omega_2$\,.

Note that in the case when we have 
a contribution from both operators $\Omega_1$ and $\Omega_2$, 
we also use the equality \ $2^{p-1}=2^p - 2^{p-1}$\,.
\qed
\enddemo

\medskip

\noindent
{\bf Example 13.} \ (i) For $n=5$, $k=1$, we have
$$
\P_{(3,2)}\, Y_{[1,3]}\, \nabla_1^D(5) = \P_2
$$
and this comes from the contribution of both operators $\Omega_1$, $\Omega_2$:
$$
\Q_{(3,2)} \, Y_{[1,3]} \, x_1 x_2 \nabla_2^B(5)=
\Q_{(3,2)} \, Y_{[2,4]} \, \nabla_2^B(5)= 4 \Q_2
$$
by Theorem 9, and
$$
\Q_{(3,2)} \, Y_{[1,3]} \, x_1 \nabla_1^B(5) \, \d_1 \d_2 \d_3 
= \Q_{(3,2)} \, Y_{[2]} \, \nabla_1^B(5)= -2 \Q_2\,.
$$
by Fact 4.
\smallskip

\noindent
(ii) \ For $n=7$, $k=2$, we have
$$
\P_{(6,5,4,3,2,1,0)}\, Y_{[0,1,2,2]} \, \nabla_2^D(7) = - \P_{(5,3)}
$$
and only the operator $\Omega_1$ gives the contribution:
$$
\Q_{(6,5,4,3,2,1)}\, Y_{[0,1,2,2]}\, x_1 x_2 x_3 x_4 \nabla_{4}^B(7) 
= \Q_{(6,5,4,3,2,1)}\, Y_{[1,2,3,3]}\, \nabla_{4}^B(7) 
=-16 \Q_{(5,3)}
$$
by Fact 4.

\smallskip

\noindent
(iii) \ For $n=7$, $k=2$, we have
$$
\P_{(6,5,4,3,2,1,0)}\, Y_{[1,1,2,3]} \, \nabla_2^D(7) = -\P_{(6,3,1)}
$$
and the contribution comes from both operators $\Omega_1$, $\Omega_2$:
$$
\Q_{(6,5,4,3,2,1)}\, Y_{[1,1,2,3]}\, x_1 x_2 x_3 x_4 \, \nabla_{4}^B(7) 
=\Q_{(6,5,4,3,2,1)}\, Y_{[2,2,3,4]}\, \nabla_{4}^B(7) 
=-16 \Q_{(6,3,1)}
$$
by Theorem 9, and
$$
\aligned
\Q_{(6,5,4,3,2,1)}\, Y_{[1,1,2,3]}\,& x_1 x_2 x_3\, \nabla_{3}^B(7)\, 
\d_1 \d_2 \d_3 \cr
&=\Q_{(6,5,4,3,2,1)}\,Y_{[2,2,3]}\, \nabla_{3}^B(7)\,
=8 \Q_{(6,3,1)}\,. \cr
\endaligned
$$
by Fact 4.

\bigskip\medskip

\head 
{\bf 4. Applications to $\P$-polynomials and orthogonal Schubert polynomials}
\endhead

\smallskip

The following presentation of a $\P$-polynomial in the form
$$
\P_I = x^{\alpha(I)}\, \P_{\rho(n-1)}\, \Omega(I)\,,
$$
where $\alpha(I) \subseteq \rho$ and $\Omega(I)$ is a divided difference 
operator, appears to be quite useful:

\proclaim{Lemma 14} \ Let $I=(i_1,\ldots, i_{\ell}>0)\subseteq 
\rho(n-1)$ be a strict partition.
If $n$ and $\ell$ are of the same parity, we set $h:=n-\ell$, and
$\{ j_1< \cdots < j_h \}: = \{ 1,\ldots,n \} \smallsetminus \{ i_1+1,\ldots,
i_{\ell}+1 \}$. 
If $n$ and $\ell$ are of different parity, we set $h:=n-\ell-1$, and
$\{ j_1< \cdots < j_h \}: = \{ 1,\ldots,n \} \smallsetminus \{ i_1+1,\ldots,
i_{\ell}+1,1 \}$.

\smallskip

Then for $\alpha(I): = [n-j_1,\ldots,n-j_h,0,\ldots,0]$ and $k:=h/2$,
$$
x^{\alpha(I)} \, \P_{\rho(n-1)}\,\partial_{[2k,\ldots,1]} \nabla_k^D(n)
= (-1)^s\, \P_I\,,
\eqno(4.1)
$$
where $s$ is the number of positions of the parts erased in 
$\rho(n-1)$ in order to get the partition $I$.
\endproclaim

The assertion of the lemma is a direct consequence of Theorem 11 and
the definition of a Schur $S$-polynomial via the Jacobi symmetrizer.

\medskip

Now, with every strict partition $I=(i_1,\ldots,i_{\ell}>0)$\,,
we associate the following element $v(I)\in \frak D_n$\,. 
If $n-\ell$ is even, we set
$$
v(I):=[i_1+1, i_2+1, \ldots, i_{\ell}+1, 
\overline{j_1}, \ldots, \overline{j_h}]\,,
\eqno(4.2)
$$
and if $n-\ell$ is odd,
$$
v(I):=[i_1+1, i_2+1, \ldots, i_{\ell}+1,1, 
\overline{j_1}, \ldots ,\overline{j_h}]\,.
\eqno(4.3)
$$
(The notation is the same as in Lemma 14.)

\proclaim{Theorem 15} \ For a strict partition $I\subseteq \rho(n-1)$, 
$$
x^{\rho} \P_{\rho(n-1)} \, \partial_{v(I)} = (-1)^{|I|+{n\choose 2}} \P_I.
\eqno(4.4)
$$
\endproclaim

The proof of this result is analogous to the proof of [LP1, Thm. A.6]. Using
the notation of Lemma 14, we have
$$
\d_{v(I)} = \d_{\sigma}\, \d_{[k,k-1,\ldots,1]}\, \nabla_k^D(n)\,,
\eqno(4.5)
$$
where
$$
\sigma=
\left\{ \matrix  [j_1,\ldots,j_h,i_1+1,\ldots,i_{\ell}+1], & n & \hbox{and} &
\ell & \hbox{of the same parity} \cr
 [j_1,\ldots,j_h,i_1+1,\ldots,i_{\ell}+1,1], & n & \hbox{and} &
\ell & \hbox{of different parity}\,. \cr
\endmatrix 
\right.
$$
Note that $x^{\rho}\, \d_{\sigma}= x^{\alpha(I)}$, and hence the assertion
follows by Lemma 14.

\smallskip

This result leads to the following characterization of $\P$-polynomials
via orthogonal divided differences:

\proclaim{Corollary 16} For a strict partition $I\subseteq \rho(n-1)$,
we set $w(I):= v(I)^{-1} w_0^D$. More explicitly,
for even $\ell$, $w(I)=[\overline{i_1+1},\ldots,\overline{i_{\ell}+1},
j_1,\ldots,j_h]^{-1}$, and for odd $\ell$,
$w(I)=[\overline{i_1+1},\ldots,\overline{i_{\ell}+1},\overline 1,j_1,
\ldots,j_h]^{-1}$. 
Then 
$w=w(I)$ is the unique element of $\frak D_n$ such that $\ell(w)=|I|$
and $\P_I \, \partial_w \ne 0$. In fact, \ $\P_I \, 
\partial_{w(I)}= (-1)^{|I|}$.

\endproclaim

This can be also seen by geometric considerations
(see [P1] and [LP1]), with the help of the characteristic map ([B], [D1,D2]).

\medskip

More generally, consider, for any $w\in \frak D_n$, 
the {\it orthogonal Schubert polynomial} 
$$
X_w^D=X_w^D(n)=
x^{\rho} \P_{\rho(n-1)} \partial_{w_0^D w}
\eqno(4.6)
$$
of degree $\ell(w)$. Arguing in the same way as in [LP1, pp.33--36],
one shows that these Schubert polynomials have the {\it stability property}
in the sense 
that for $w\in \frak D_n \subset \frak D_{n+1}$,
$$
X_w^D(n+1)|_{x_{n+1}=0} = X_w^D(n).
\eqno(4.7)
$$
Together with the ``maximal Grassmannian property"
from Theorem 15, asserting that, for even $\ell$,
$$
X_{[\overline{i_1+1},\ldots,\overline{i_{\ell}+1},j_1,\ldots,j_h]}^D=
(-1)^{|I|+{{n}\choose 2}} \P_I\,,
\eqno(4.8)
$$
and for odd $\ell$,
$$
X_{[\overline{i_1+1},\ldots,\overline{i_{\ell}+1},\overline 1,
j_1,\ldots,j_h]}^D=
(-1)^{|I|+{{n}\choose 2}} \P_I\,,
\eqno(4.9)
$$
this shows that they provide a natural tool for the
cohomological study of Schubert varieties for the orthogonal group
$SO(2n)$ and the related degeneracy loci.

\medskip

We also record  

\proclaim{Proposition 17} \ For a strict partition $I=(i_1,i_2,i_3,i_4, \ldots)
\subseteq \rho(n-1)$,
$$
\P_I \, \dd \d_2 \cdots \d_{i_1}
\d_1 \d_2 \cdots \d_{i_2} = (-1)^{i_1+i_2} \P_{(i_3, i_4, \ldots)}.
\eqno(4.10)
$$
\endproclaim

To see this, we argue in a manner similar to the proof of [LP1, Prop. 5.12].
For $J=(i_3,i_4,\ldots)$, we choose the presentation from Lemma 14:
$$
\pm \P_I = x^{\alpha(I)}\, \P_{\rho(n-1)} \d_u \ \ , \ \ \ \ 
\pm \P_J = x^{\alpha(J)}\, \P_{\rho(n-1)} \d_v 
$$
for appropriate $u, v \in \frak D_n$. Let $\sigma \in \frak S_n$ be the
permutation such that 
$$
v=\sigma\, u\, \ss s_2 \cdots s_{i_1}
s_1 s_2 \cdots s_{i_2}\,.
$$ 
The assertion now follows from \ $x^{\alpha(I)} \d_{\sigma}= x^{\alpha(J)}$.

\bigskip\medskip

\head
{\bf Appendix: results in type B}
\endhead

In this appendix we give a summary of the results for type $B_n$.
They are obtained directly from the results for type $C_n$ in [LP1], by changing
$\d_0$ in loc.cit. to $-2 \d_0$, and read as follows:
writing $\nabla:= \nabla_n^B(n)$, we have

\proclaim {Theorem 18} \ (i) \ For $\alpha\in \N^n,\, \alpha \subseteq \rho$,
$$
\ Y_{\alpha}\, \P_{\rho(n)} \, \nabla = 
(-1)^{|\alpha|+{{n+1}\choose 2}} Y_{\alpha}^{\omega}.
\eqno(5.1)
$$
\medskip

\noindent 
(ii) \ For strict $I\varsubsetneq \rho(n)$, and $\alpha \subseteq \rho$,
$$
\ Y_{\alpha} \, \P_I \, \nabla = 0 \ .
\eqno(5.2)
$$
\endproclaim

\bigskip

Denote by
$$
\langle \ ,\ \rangle : \Sym(n) \times \Sym(n) \to Sym^D(n) 
$$
the scalar product defined for $f, g\in Sym(n)$ by
$$
\langle f,g \rangle := f g \, \nabla\,.
\eqno(5.3)
$$
For strict partitions $I,J\subseteq \rho(n)$, one has
$$
\langle \P_I,\P_{\rho(n) \smallsetminus J} \rangle 
= (-1)^{{{n+1}\choose 2}}\, \delta_{I J}, 
\eqno(5.4)
$$
where $\rho(n) \smallsetminus J$ is the strict partition whose parts 
complement the parts of $J$ in $\{n,n-1,\ldots, 1\}$ \ (cf. [PR]).

Consequently, the polynomial ring $\Pol= \Bbb Z [{1 \over 2}][x_1,\ldots, x_n]$
is a free $\Sym^B(n)$-module with basis $Y_{\alpha}\, \P_I$, where
$\alpha$ ranges over subsequences contained in $\rho$ and $I$
runs over all strict partitions contained in $\rho(n)$. 
Note that the element of the maximal degree of this basis is
$x^{\rho} \P_{\rho(n)}$. 
Let
$$
[ \ , \ ]: \Pol \times \Pol \to \Sym^B(n)
$$
be a scalar product, defined for $f,g \in \Pol$ by
$$
[f,g]:= f g \,\d_{w_0^B}.
\eqno(5.5)
$$
One has, for $\alpha, \beta \subseteq \rho$ and strict partitions $I,J\subset
\rho(n)$,
$$
\bigl[Y_{\alpha}^{\omega}\, \P_I , 
Y_{\beta'}\, \P_{\rho(n)\smallsetminus J} \bigr] = 
(-1)^{|\alpha| + {{n+1}\choose 2}}\, \delta_{\alpha \beta}\, \delta_{I J}\,.
\eqno(5.6)
$$ 
(See (2.5).)

\medskip
  
Let $\Y=\{ y_1,\ldots,y_n\}$ be a second set of indeterminates of
cardinality $n$. 
The symbol $\equiv$ will mean: \ ``congruent modulo the ideal
generated by the relations \hfill\break
$ f(x_1^2,\ldots, x_n^2) =
f(y_1^2,\ldots, y_n^2)$, where $f\in Sym(n)$."

\medskip

Following Fulton [F2,F3], define 
$$
\aligned
F(\X,\Y) &:=\ |\P_{n+1+j-2i}(\X)+\P_{n+1+j-2i}(\Y)|_{1\le i,j\le n} \\
& \\
&=\ {\left|\ \CD
\P_{n}(\X)+\P_{n}(\Y) & 0 &.\;.\;. & \\
\P_{n-2}(\X)+\P_{n-2}(\Y) &\ \P_{n-1}(\X)+\P_{n-1}(\Y) &\ .\;.\;.\ & \\
\vbox{\offinterlineskip\hbox{.}\vskip3pt\hbox{.}\vskip3pt\hbox{.}}
&\vbox{\offinterlineskip\hbox{.}\vskip3pt\hbox{.}\vskip3pt\hbox{.}}
&\vbox{\offinterlineskip\hbox{.}\vskip3pt\hbox{\ \ .}\vskip3pt
\hbox{\ \ \ \ .}} & \\
&& 1 &&  \P_1(\X)+\P_1(\Y) \\
\endCD\ \right|}\,. \\
\endaligned
\eqno(5.7)
$$

Following [PR], define
$$ 
\P(\X, \Y):= \sum \P_I(\X)\, \P_{\rho(n) \smallsetminus I}(\Y) \ ,  
\eqno(5.8)
$$
where the summation is over all strict partitions $I\subseteq \rho(n)$.
The reasoning in [LP1, Sect.2] made for case $C_n$ adapts to case $B_n$ 
and furnishes:

\proclaim {Proposition 19} \ We have 

\noindent
(i) \ $$F(\X,\Y) \equiv \P(\X,\Y). \eqno(5.9)$$ 

\noindent
(ii) \ For every $w\in \frak B_n\smallsetminus \S_n$,
$$
\P(\X^w, \X)=0\,,
\eqno(5.10)
$$
and for every $w\in \S_n$,
$$
\P(\X^w, \X)= \P(\X, \X)= s_{\rho(n)}(\X)\,.
\eqno(5.11)
$$

\noindent
(iii) \ For every $f\in \Sym(n)$, 
$$
\langle f(\X)\,, F(\X,\Y) \rangle \equiv (-1)^{{{n+1} \choose 2}}\, f(\Y) \ . 
\eqno(5.12)
$$

\noindent
(iv) \ For every $f\in \Pol$,
$$
\bigl[ f(\X)\,, \prod_{n \ge i>j \ge 1}(x_i-y_j)\, F(\X,\Y) \,\bigr] 
\equiv f(\Y)\,.
\eqno(5.13)
$$
\endproclaim
In other words, $F(\X,\Y)$ is a reproducing kernel for the scalar product
$\langle \ , \ \rangle$, and $\prod_{i>j}(x_i-y_j)\, F(\X,\Y)$ is a reproducing
kernel for $[ \ , \ ]$. One can show that the ``vanishing property"
{\it (ii)} characterizes $\P(\X,\Y)$ up to $\equiv$. The congruence
{\it (i)} can be also derived from geometry by comparing the classes
of diagonals in flag bundles associated with $SO(2n+1)$ given in
[F2,F3] and [PR] (see also [G]).

\medskip

\proclaim{Proposition 20}
Suppose $n\ge p>0$. Let  $IpJ\subseteq \rho(n)$ be a strict partition 
and $H\subseteq \rho(n)$ a strict partition not containing $p$.
Then
$$
x_1^{n-p}\P_{IpJ} \, \d_0 \d_1 \cdots \d_{n-1}
=(-1)^{\ell(I)+n} \P_{IJ} 
\eqno(5.14)
$$
and
$$
x_1^{n-p}\P_H \, \d_0 \d_1 \cdots \d_{n-1} = 0 \,.
\eqno(5.15)
$$
\endproclaim

More generally,

\proclaim{Theorem 21} \quad Let $0<k\leq n$ and let 
$\alpha=[\alpha_1\leq \cdots \leq \alpha_k]\in \N^k$ 
be such that $\alpha_k \le n-k$.
Suppose that $I\subseteq \rho(n)$ is a strict partition. 
Then the image of $ \P_I \, Y_\alpha$ 
under $\nabla_k^B(n)$ is 0 
unless  $n-\alpha_1-0,\ldots, n-\alpha_k-(k-1)$ are
parts of $I$. In this case, the image is $(-1)^{k(n-1)+s} \P_J$, where 
$J$ is the strict partition with parts
$$\{i_1,\ldots, i_{\ell(I)} \} \smallsetminus
\{ n-\alpha_1-0,\ldots, n-\alpha_k-(k-1) \},
$$
and $s$ is the sum of positions of the parts erased in $I$.  
\endproclaim 

(This is a restatement of [LP1, Prop. 5.9].)

\medskip

\proclaim{Proposition 22} \quad For a strict partition $I=(i_1,i_2,\ldots)$,
$$
\P_I\, \d_0 \d_1 \cdots \d_{i_1-1} = (-1)^{i_1} \P_{(i_2,\ldots)}\,.
\eqno(5.16)
$$
\endproclaim

Now, let us associate with every strict partition $I=(i_1,\ldots,i_{\ell}>0)$
the following element of $\frak B_n$:
$$
v(I):=[i_1,\ldots,i_{\ell},\overline{j_1}, \ldots , \overline{j_h}]\,,
\eqno(5.17)
$$
where $j_1 < \cdots < j_h$\,.

\smallskip

\proclaim{Theorem 23} \ For every strict partition $I\subseteq \rho(n)$,
$$
x^{\rho} \P_{\rho(n)} \, \partial_{v(I)} = (-1)^{|I|+{{n+1}\choose 2}} \P_I.
\eqno(5.18)
$$
\endproclaim

This leads to the following characterization of $\P$-polynomials
via divided differences:

\proclaim{Corollary 24} \ 
For any strict partition $I$, let $w(I):= v(I)^{-1} w_0^B$, that is
$w(I)= [\overline{i_1},\ldots,\overline{i_{\ell}},j_1,\ldots,j_h]^{-1}$\,.
Then $w=w(I)$
is the unique element of $\frak B_n$ such that $\ell(w)=|I|$ and
$\P_I \, \partial_w \ne 0$. In fact, $\P_I \, \partial_{w(I)}=(-1)^{|I|}$. 
\endproclaim

This can be also seen by geometric considerations (see [P1] and [LP1]),
with the help of the characteristic map ([B], [D1,D2]).

\smallskip

More generally, consider, for any $w\in \frak B_n$, the 
{\it orthogonal Schubert polynomial} 
$$
X_w^B= X_w^B(n) = x^{\rho} \P_{\rho(n)} \partial_{w_0^B w}
\eqno(5.19)
$$
of degree $\ell(w)$. Arguing in the same way as in [LP1, pp.33--36],
one shows that these Schubert polynomials have the {\it stability property}
in the sense that for $w\in \frak B_n \subset \frak B_{n+1}$,
$$
X_w^B(n+1)|_{x_{n+1}=0} = X_w^B(n).
\eqno(5.20)
$$
Together with the ``maximal Grassmannian property"
from Theorem 23, asserting that
$$
X_{[\overline{i_1},\ldots,\overline{i_{\ell}},j_1,\ldots,j_h]}^B=
(-1)^{|I|+{{n+1}\choose 2}} \P_I\,,
\eqno(5.21)
$$
this shows that they provide a natural tool for the
cohomological study of Schubert varieties for the orthogonal group
$SO(2n+1)$ and the related degeneracy loci.

\bigskip

We now give

\smallskip

\noindent
{\bf Proof of Theorem 9}

\smallskip

Given a symmetric function $f$, let $D_f$ be the Foulkes derivative i.e.
the adjoint 
operator to the multiplication by $f$ w.r.t. the standard scalar product
on the ring $\Sym$ of symmetric functions in a countable number of
variables (cf. [M1]). 
We use the following vertex operators on $\Sym$:
$$
U^s := 1 - D_{P_1} s_1 + D_{P_2} s_2 - \cdots \,,
\eqno(5.22)
$$
$$
U^e := 1 - D_{P_1} e_1 + D_{P_2} e_2 - \cdots \,,
\eqno(5.23)
$$
and
$$
V^e:=1 - D_{e_1} P_1+ D_{e_2} P_2 - \cdots \,.
\eqno(5.24)
$$
We refer to [LP1, p.24] for the definitions of Schur $P$-functions $P_I$
[S]. In loc. cit. the reader can also find a definition of 
$Q'$-functions $Q'_I$ [LLT2], used in the following proposition:    

\proclaim {\bf Proposition 25} \ Let $I$ be a strict partition.
We have the following identities 
of symmetric functions in $\Sym$:
$$
\Q_I \, U^s
= \left\{ \matrix \Q_I, & \ell(I) & \hbox{even} \cr
                0, &\ell(I) &\hbox{odd}\,; \cr
\endmatrix 
\right.
\eqno(5.25)
$$
$$
Q_I' \, U^e
= \left\{ \matrix Q_I', & \ell(I) & \hbox{even} \cr
                0, &\ell(I) &\hbox{odd}\,; \cr
\endmatrix 
\right.
\eqno(5.26)
$$
and
$$
P_I \, V^e
= \left\{ \matrix P_I, & \ell(I) & \hbox{even} \cr
                0, &\ell(I) &\hbox{odd}\,. \cr
\endmatrix 
\right.
\eqno(5.27)
$$
\endproclaim
 
\demo{Proof}
First of all,
arguing as in [LP1, pp.24--27], with the help of the operators $V^e$,
$U^s$, and $U^e$ instead of $V^e_k$, $U^s_k$, and $U^e_k$ in loc. cit., 
we note that the equalities (5.25), (5.26), and (5.27) are equivalent.  

\smallskip

We show (5.27). It suffices to prove the statement when the set of 
indeterminates $\{x_1,\ldots, x_n\}$ is of finite cardinal
$n>|I|$\,. 

Besides the well-known equality: for $k>0$,
$$
x_1^k \prod_{2\le i \le n} (x_1+x_i)\, \d_1 \d_2 \cdots \d_{n-1}
= P_k(x_1,\ldots,x_n)\,,
\eqno(5.28)
$$
we also need the following formula from [P2]:
\proclaim {Fact 26} \ For a strict partition $I$,
$$\aligned
P_I(x_2,\ldots,x_n)& \, \prod_{2\le i \le n}(x_1 + x_i)\, \d_1 \d_2 
\cdots \d_{n-1} \cr
&= \left\{ \matrix (-1)^{n-1} P_I(x_1,\ldots,x_n), & n-\ell(I) 
& \hbox{odd} \cr
                0, &n-\ell(I) &\hbox{even}. \cr
\endmatrix 
\right.
\endaligned
\eqno(5.29)
$$
\endproclaim
\noindent
(More precisely, Eq.(5.29) is a special case of the following formula given in 
[P2, Prop.1.3(ii)]: \ 
Let $q$, $r$, $k$, and $h$ be integers
such that $0<q<n$, $n=q+r$, $0 \le k \le q$, and $0 \le h \le r$.
Suppose $I=(i_1,\ldots,i_k)\in \Bbb N^{* k}$ and $J=(j_1,\ldots,j_h)\in 
\Bbb N^{* h}$. Then
$$
P_I(x_1,\ldots,x_q)\, P_J(x_{q+1},\ldots,x_n)
\prod_{1\le i \le q<j\le n} (x_i + x_j)\,
 \matrix 
\d_q & \cdots & \d_{n-1} \cr
\vdots & & \vdots \cr
\d_2 & \cdots & \d_{r+1} \cr
\d_1 & \cdots & \d_r \cr
\endmatrix 
$$
$$
= d \cdot P_{(i_1,\ldots,i_k,j_1,\ldots,j_h)}(x_1,\ldots,x_n)\,, 
\eqno(5.30)
$$
where $d$ is zero if $(q-k)(r-h)$ is odd and
$$
d=(-1)^{(q-k)r} {{\lfloor (n-k-h)/2 \rfloor}\choose{\lfloor (q-k)/2 \rfloor}}
\eqno(5.31)
$$
otherwise.

\smallskip

We get Eq.(5.29) as Eq.(5.30) specialized to $q=1$ and $k=0$.)

\bigskip

To end the proof of (5.27), we first write
$$
P_I(x_2,\ldots,x_n)=P_I - P_I D_{e_1} \cdot x_1 + P_I D_{e_2} \cdot x_1^2 -
P_I D_{e_3} \cdot x_1^3 + \cdots \,,
\eqno(5.32)
$$
where the RHS is evaluated in the first $n$ variables. Then we multiply
both sides of (5.32) by
$$
(x_1+x_2)(x_1+x_3) \cdots (x_1+x_n)
$$
and apply the operator \ $\d_1 \d_2 \cdots \d_{n-1}$. We get the following
equalities of symmetric polynomials in the first $n$ variables.
If $n$ is odd, the RHS of the so-obtained equality becomes  
$$
P_I - P_I D_{e_1}\cdot P_1 + P_I D_{e_2}\cdot P_2 - P_I D_{e_3}\cdot P_3 
+ \cdots
$$
by (5.28) and (5.29), and its LHS is equal to
$$
\left\{ \matrix  P_I,  & \ell(I) & \hbox{even} \cr
                0, & \ell(I) &\hbox{odd} \cr
\endmatrix 
\right.
$$
by (5.29). This shows (5.27) for odd $n$.
If $n$ is even, the RHS of the obtained equality becomes
$$
0 - P_I D_{e_1}\cdot P_1 + P_I D_{e_2}\cdot P_2 - P_I D_{e_3}\cdot P_3 + \cdots
$$
by (5.28) and (5.29), and its LHS is equal to
$$
\left\{ \matrix  0,  & \ell(I) & \hbox{even} \cr
                -P_I, & \ell(I) &\hbox{odd} \cr
\endmatrix 
\right.
$$
by (5.29). This shows (5.27) for even $n$. 
\smallskip

Thus the proposition has been proved.
\qed
\enddemo

\medskip

Let now $\d_0^C$ be the divided difference defined by
$$ 
\Pol \ni f \mapsto f\, \d_0^C :=  (f-f^{s_0})/{(2x_1)}\ .   
\eqno(5.33)
$$
\smallskip

\noindent
Arguing similarly as in [LP1, pp.27--28], using 
$$
\d_0^C = D_{P_1} - D_{P_2}\, x_1 + D_{P_3}\, x_1^2 - \cdots \,,
\eqno(5.34)
$$
and the formula 
$$
x_1^p \, \d_1 \cdots \d_{n-1} = s_{p-n+1} (x_1,\ldots, x_n)\,,
\eqno(5.35)
$$
one shows

\proclaim {Lemma 27} \ As operators on $\Sym$, evaluated in
symmetric polynomials in the first $n$ variables, 
$$
1 - U^s = \d_0^C\, x_1^n\, \d_1 \cdots \d_{n-1}\,.
\eqno(5.36)
$$
\endproclaim

\smallskip

Equations (5.25) and (5.36), together with the equalities
$$
x_1^{2m}\, \d_0^C = \d_0^C \, x_1^{2m} \ \ \ \hbox{and} \ \ \
x_1^{2m+1}\, \d_0^C = -\d_0^C \, x_1^{2m+1} + x_1^{2m}\,,
\eqno(5.37)
$$
applied for even $n=2m$ or odd $n=2m+1$ accordingly, imply 

\proclaim {Proposition 28} \ Let $I\subseteq \rho(n)$ be a strict
partition. We have 
$$
\Q_I \, x_1^n \, \d_0 \d_1 \cdots \d_{n-1}
= \left\{ \matrix -2\Q_I, & n+\ell(I) & \hbox{odd} \cr
                0, &n+\ell(I) &\hbox{even}\,. \cr
\endmatrix 
\right.
\eqno(5.38)
$$
\endproclaim

\smallskip

Eq.(5.38) is the content of Theorem 9 for $k=1$. For higher $k$, one gets
the desired assertion by [LP1, Thm 5.1], Proposition 28, and [LP1, Lemma 5.10].
(Note that this last fact holds true for any nonnegative integer $\alpha_1$,
in the notation of loc. cit., as is clear from its proof.)

\smallskip

This ends the proof of Theorem 9.

\bigskip

Finally, we take this opportunity to correct some misprints in [LP1]:
 -- should read:

\smallskip
\noindent
p.$11_{11}$ \  ``... a partition ..."

\vskip5pt
\noindent
p.$13_2$ \ `` $\langle \ , \ \rangle : \Cal S \Cal P(X) \times \Cal S \Cal P(X)
\to \Cal S \Cal P(X^2)$ "

\vskip5pt
\noindent
p.$36_2$ \ ``... $\partial_u' \Cal C_w = \Cal C_v$ ..."

\vskip5pt
\noindent
p.$37^4$ \ ``... is $\nabla_k \circ \partial_{\omega^{(k)}}'$ ..."

\vskip5pt
\noindent
p.$37_{10}$ \ ``... $\partial_{w_I}'(\Q_I(X))=1$ ... "

\vskip5pt
\noindent
Moreover, in Example 5.11, the sequence of successive signs is: +, +, --, 
-- \,.

{\eightrm

\bigskip

\noindent
A.L. : C.N.R.S., Institut Gaspard Monge, Universit\'e de Marne-la-Vall\'ee,
5 Bd Descartes, Champs sur Marne, 77454 Marne La Vall\'ee Cedex 2, France

\noindent
e-mail: Alain.Lascoux@univ-mlv.fr

\bigskip

\noindent
P.P. : Mathematical Institute of Polish Academy of Sciences, Chopina 12,
87-100 Toru\'n, Poland

\noindent
e-mail: pragacz@mat.uni.torun.pl

\enddocument
\bye